\documentclass{amsart}
\usepackage{amsmath,amssymb,latexsym,mathrsfs}
\usepackage[mathscr]{eucal}
\usepackage{enumerate}
\usepackage{color}

\theoremstyle{remark}

\input{tcilatex}

\begin{document}
\title{Monotonic Extensions of Lipschitz Maps}
\author{Efe A. Ok}
\address{Department of Economics and Courant Institute of Mathematical
Sciences, New York University}
\email{efe.ok@nyu.edu}
\date{April 18, 2023}

\begin{abstract}
We study the problem of extending an order-preserving real-valued Lipschitz
map defined on a subset of a partially ordered metric space without
increasing its Lipschitz constant and preserving its monotonicity. We show
that a certain type of relation between the metric and order of the space,
which we call \textit{radiality}, is necessary and sufficient for such an
extension to exist. Radiality is automatically satisfied by the equality
relation, so the classical McShane-Whitney extension theorem is a special
case of our main characterization result. As applications, we obtain a
similar generalization of McShane's uniformly continuous extension theorem,
along with some functional representation results for radial partial orders.
\end{abstract}

\maketitle

\section{Introduction}

The most important continuous extension theorem for order-preserving
functions is \textit{Nachbin's extension theorem}. This theorem considers a
partially ordered topological space, and gives conditions under which a
continuous and increasing real-valued function defined on a compact subset
of such a space can be extended to the entire space in such a way to remain
continuous and increasing. It has found profound applications, especially in
the field of decision theory. (The references for the present discussion are
provided in the body of the paper.)

Another extension theorem of great importance is the famous \textit{%
McShane-Whitney extension theorem}. This theorem shows that any Lipschitz
map defined on a subset of a metric space can be extended to the entire
space without increasing the Lipschitz constant of the original map. This
theorem paved the way toward various types of Lipschitz extension theorems
for Banach space-valued Lipschitz maps, presently a topic of active research
in geometric functional analysis. In addition, it has recently been
pivotally used in the literature on machine learning and metric data
analysis.

The primary objective of this note is to understand to what extent a Nachbin
type generalization of the McShane-Whitney theorem is possible. To state our
query precisely, consider a 1-Lipschitz real-valued function $f$ on a subset 
$S$ of a metric space $X.$ Now suppose $X$ is endowed with a partial order $%
\succcurlyeq ,$ and that $f$ is order-preserving (in the sense that $%
f(x)\geq f(y)$ for every $x,y\in S$ with $x\succcurlyeq y$). The problem is
to determine under what conditions (that do not depend on $S$ and $f$), one
can extend $f$ to an order-preserving 1-Lipschitz map on $X.$ Our main
result (Theorem 2, below) says that this is possible if, and only if, $%
\succcurlyeq $ satisfies a rather demanding condition, which we call \textit{%
radiality.} (We actually prove a slightly more general result that covers $%
\ell _{\infty }(T)$ valued-maps as well, for any nonempty set $T$.)
Radiality of $\succcurlyeq $ demands that if $x\succcurlyeq y$ while $%
z\succcurlyeq y$ does not hold, then the distance between $x$ and $z$ is
larger than that between $y$ and $z$ (and similarly for the case where not $%
y\succcurlyeq x$ but $y\succcurlyeq z$). While it is obviously strong, this
condition is necessary for the sought monotonic Lipschitz extension theorem.
Moreover, when $\succcurlyeq $ is total, it reduces to \textit{radial
convexity} which is commnoly used in the field of topological order theory.
Finally, the equality ordering is radial, so our extension theorem
generalizes the McShane-Whitney extension theorem (just like Nachbin's
theorem generalizes the Tietze extension theorem).

We also present some applications of our monotonic Lipschitz extension
theorem. First, we show that every radial partial order on a (compact)
metric space can be represented by means of a (compact) collection of
order-preserving Lipschitz functions. An immediate corollary of this is that
every radial partial order is closed. Second, we prove that on any radial
partially ordered $\sigma $-compact metric space $X$, there is a strictly
increasing Lipschitz map $F$ (in the sense that $F(x)>F(y)$ for every
distinct $x,y\in X$ with $x\succcurlyeq y$). Finally, we combine our
extension theorem with the recent remetrization approach introduced in Beer 
\cite{BeerBook} to show that if $\succcurlyeq $ is a radial partial order on
a metric space $X,$ then any bounded (or more generally, Lipschitz for large
distances), order-preserving and uniformly continuous map on a subset of $X$
can be extended to a function on the entire space in such a way that it
remains order-preserving and uniformly continuous. Radiality can actually be
relaxed substantially in this result, but characterizing those metric posets
on which such an extension is possible is presently an open problem.

\section{Preliminaries}

\subsection{Posets}

Let $X$ be a nonempty set. A \textit{preorder} on $X$ is a reflexive and
transitive binary relation on $X,$ while a \textit{partial order} on $X$ is
an antisymmetric preorder on $X.$ We refer to the ordered pair $%
(X,\succcurlyeq )$ as a \textit{poset }if $\succcurlyeq $ is a partial order
on $X.$ (In this context, $X$ is called the \textit{carrier }of the poset.)
A preorder on $X$ is \textit{total} if any two elements $x$ and $y\,$of $X$
are $\succcurlyeq $\textit{-comparable}, that is, either $x\succcurlyeq y$
or $y\succcurlyeq x$ holds. A total partial order on $X$ is said to be a 
\textit{linear order} on $X;$ in this case, we refer to $(X,\succcurlyeq )$
as a \textit{loset. }

Let $(X,\succcurlyeq )$ be a poset. For any $x\in X,$ we define $%
x^{\downarrow }:=\{z\in X:x\succcurlyeq z\}$ and $x^{\uparrow }:=\{z\in
X:z\succcurlyeq x\}$. (A set of the former type is said to be a \textit{%
principle ideal }in $(X,\succcurlyeq ),$ and one of the latter type is
called a \textit{principal filter} in $(X,\succcurlyeq )$.) In turn, for any 
$S\subseteq X,$ we define the $\succcurlyeq $\textit{-decreasing closure} of 
$S$ as $S^{\downarrow }:=\tbigcup\nolimits_{x\in S}x^{\downarrow }$, and
define the $\succcurlyeq $\textit{-increasing closure} $S^{\uparrow }$ of $S$
dually. In turn, $S$ is said to be $\succcurlyeq $\textit{-decreasing }if $%
S=S^{\downarrow },$ and $\succcurlyeq $\textit{-increasing }if $%
S=S^{\uparrow }$.

Given any poset $(X,\succcurlyeq ),$ we denote the asymmetric part of $%
\succcurlyeq $ by $\succ $, that is, $x\succ y$ means $y\neq x\succcurlyeq y.
$ We also define the binary relation $\succcurlyeq ^{\bullet }$ on $X$ as%
\begin{equation*}
x\succcurlyeq ^{\bullet }y\hspace{0.2in}\text{iff\hspace{0.2in}not }%
y\succcurlyeq x.
\end{equation*}%
Thus $x\succcurlyeq ^{\bullet }y$ means that either $x\succ y$, or $x$ and $y
$ are not $\succcurlyeq $-comparable. It is plain that $\succcurlyeq
^{\bullet }$ is an irreflexive relation. In general, this relation is
neither symmetric nor asymmetric, nor it is transitive. When $\succcurlyeq $
is total, however, $\succcurlyeq ^{\bullet }$ equals $\succ $.

A function $f:X\rightarrow Y$ from a poset $X=(X,\succcurlyeq )$ to a poset $%
Y=(Y,\trianglerighteq )$ is said to be \textit{order-preserving} if for
every $x,y\in X$, 
\begin{equation*}
x\succcurlyeq y\hspace{0.2in}\text{implies\hspace{0.2in}}f(x)%
\trianglerighteq f(y).
\end{equation*}%
If $Y$ is $(\mathbb{R},\geq ),$ where $\geq $ is the usual order, we refer
to $f$ simply as $\succcurlyeq $\textit{-increasing. }Note that the
indicator function of any $\succcurlyeq $-increasing subset of $X$ is an $%
\succcurlyeq $-increasing map.

\subsection{Normally Ordered Topological Posets}

A \textit{topological poset }is an ordered pair $(X,\succcurlyeq )$ where $X$
is a topological space and $\succcurlyeq $ is a partial order on $X$ such
that $\succcurlyeq $ is closed in $X\times X$ (relative to the product
topology). On the other hand, a \textit{normally ordered topological space}
is an ordered pair $(X,\succcurlyeq )$ where $X$ is a topological space and $%
\succcurlyeq $ is a partial order on $X$ such that for every disjoint closed
subsets $A$ and $B$ such that $A$ is $\succcurlyeq $-decreasing and $B$ is $%
\succcurlyeq $-increasing, there exist disjoint open subsets $O$ and $U$ of $%
X$ such that $O$ is $\succcurlyeq $-decreasing and contains $A$, and $U$ is $%
\succcurlyeq $-increasing and contains $B$. If, in addition, $\succcurlyeq $
is closed in $X\times X,$ we refer to $(X,\succcurlyeq )$ as a \textit{%
normally ordered topological poset.}

In his seminal work, Nachbin \cite{Na} has studied such spaces and obtained
the following generalization of the classical Tietze extension theorem:

\bigskip 

\noindent \textsf{\textbf{The Nachbin Extension Theorem.}}\textbf{\ }\textit{%
Let }$(X,\succcurlyeq )$\textit{\ be a normally ordered topological poset.} 
\textit{Then for every compact subset }$S$\textit{\ of }$X$\textit{, and
every }$\succcurlyeq $\textit{-increasing and continuous }$f:S\rightarrow 
\mathbb{R}$\textit{, there is an }$\succcurlyeq $\textit{-increasing and
continuous }$F:X\rightarrow \mathbb{R}$\textit{\ with }$F|_{S}=f.$

\bigskip 

This is a truly powerful extension theorem which holds true also when $%
\succcurlyeq $ is not antisymmetric (cf. \cite{Min}). It is used extensively
in decision theory; see, for instance, \cite{Bosi}, \cite{E-H}, \cite{E-O},
and references cited therein. 

\subsection{Partially Ordered Metric Spaces}

A \textit{partially }(resp., \textit{linearly}) \textit{ordered metric space 
}is an ordered triplet $(X,d,\succcurlyeq )$ such that $(X,d)$ is a metric
space and $(X,\succcurlyeq )$ is a poset (resp., loset). If, in addition, $%
\succcurlyeq $ is a closed subset of $X\times X,$ we refer to $%
(X,d,\succcurlyeq )$ as a \textit{metric poset} (resp., \textit{metric loset}%
).

A partially ordered metric space $(X,d,\succcurlyeq )$ is said to be \textit{%
radially convex} (or that the partial order $\succcurlyeq $ on $(X,d)$ is 
\textit{radially convex}) if 
\begin{equation*}
x\succ y\succ z\hspace{0.2in}\text{implies\hspace{0.2in}}d(x,z)\geq \max
\{d(x,y),d(y,z)\}
\end{equation*}%
for every $x,y,z\in X.$ This concept builds an appealing connection between
the order and metric structures to be imposed on a given set. Indeed, such
partially ordered metric spaces have received some attention in topological
order theory (cf. \cite{B-O}, \cite{Carruth} and \cite{Ward}, among many
others), and are often used in the topological analysis of smooth dendroids;
(cf. \cite{Fugate}).

In what follows we will need to work with a strengthening of radial
convexity. We say that a partially ordered metric space $(X,d,\succcurlyeq )$
is \textit{radial }(or that the partial order $\succcurlyeq $ on $(X,d)$ is 
\textit{radial}) if 
\begin{equation}
x\succcurlyeq ^{\bullet }y\succ z\hspace{0.2in}\text{implies\hspace{0.2in}}%
d(x,z)\geq d(x,y)  \label{d1}
\end{equation}%
and%
\begin{equation}
x\succ y\succcurlyeq ^{\bullet }z\hspace{0.2in}\text{implies\hspace{0.2in}}%
d(x,z)\geq d(y,z)\text{.}  \label{d2}
\end{equation}%
While radiality is more demanding than radial convexity, these concepts
coincide when the partial order at hand is total.

\bigskip

\noindent \textsf{\textbf{Lemma 1. }}\textit{A linearly ordered metric space
is radial if and only if it is radially convex.}

\subsection{Examples of Radial Metric Posets}

If we order and metrize any nonempty subset of $\mathbb{R}$ in the usual
way, we obtain a radial metric loset. Besides, it is plain that every
partially ordered discrete metric space is radial, and the equality relation
on any metric space is radial. But easy examples show that ordering $\mathbb{%
R}^{2}$ coordinatewise and endowing it with the Euclidean metric yields a
radially convex metric poset which is not radial.

Before proceeding further, we present a few more examples.

\bigskip 

\noindent \textit{Example 1. }Consider the poset $(X,\succcurlyeq )$ where $%
X:=\{x_{1},x_{2},x_{3},x_{4}\},$ $x_{1}\succ x_{2}\succ x_{4}$, $x_{1}\succ
x_{3}\succ x_{4}$, and $x_{2}$ and $x_{3}$ are not $\succcurlyeq $%
-comparable. (This poset is isomorphic to $(2^{S},\supseteq )$ for any
doubleton $S$.) For any $a,b\in (0,1),$ define $d_{a,b}:X\times X\rightarrow
\lbrack 0,1]$ by the matrix%
\begin{equation*}
\left[ 
\begin{array}{cccc}
0 & a & a & 1 \\ 
a & 0 & b & 1-a \\ 
a & b & 0 & 1-a \\ 
1 & 1-a & 1-a & 0%
\end{array}%
\right] 
\end{equation*}%
whose $ij$th term is $d_{ab}(x_{i},x_{j}),$ $i,j=1,...,4.$ Then, $d_{a,b}$
is a metric on $X$ iff $\min \{a,1-a\}\geq \frac{1}{2}b$. In fact, under
this parametric restriction, $(X,d_{a,b},\succcurlyeq )$ is a radially
convex metric poset. In addition, if $1-a<b<a,$ this metric poset satisfies
the condition (\ref{d2}), but not (\ref{d1}), while if $a<b<1-a,$ then the
opposite situation ensues. (In particular, this shows that there is no
redundancy in our definition of radiality.) Consequently, $%
(X,d_{a,b},\succcurlyeq )$ is a radial metric poset, provided that $\min
\{a,1-a\}\geq b.$ $\square $

\bigskip 

\noindent \textit{Example 2. }Let $T$ be a tree with a finite set $X$ of
vertices and root $x_{0}\in X.$ The \textit{path-metric} on $X$ (induced by $%
T)$ is defined as%
\begin{equation*}
\rho _{T}(x,y):=\text{the length of the path between }x\text{ and }y\text{
in }T.
\end{equation*}%
(Since $T$ is a tree, there is a unique path between any of its two
vertices.) We define $d_{T}:X\times X\rightarrow \{0,1,2\}$ by setting $%
d_{T}(x,y):=\min \{\rho _{T}(x,y),2\}$ if $x$ and $y$ are on the same path
whose one endpoint is $x_{0},$ and $d_{T}(x,y):=1$ otherwise. It is readily
checked that $d_{T}$ is a metric on $X.$ Finally, we define the partial
order $\succcurlyeq $ on $X$ by%
\begin{equation*}
x\succcurlyeq y\hspace{0.2in}\text{iff\hspace{0.2in}}y\text{ is on the path
between }x_{0}\text{ and }x.
\end{equation*}%
Then, $(X,\rho _{T},\succcurlyeq )$ is a radially convex metric poset (which
need not be radial), while $(X,d_{T},\succcurlyeq )$ is a radial metric
poset. $\square $

\bigskip

\noindent \textit{Example 3. }Let $A$ and $B$ be two disjoint bounded
subsets of a metric space $(Y,d)$. Let $\succcurlyeq _{A}$ and $\succcurlyeq
_{B}$ be radially convex linear orders on $(A,d)$ and $(B,d),$ respectively.
Let $\succcurlyeq $ be the disjoint sum of $\succcurlyeq _{A}$ and $%
\succcurlyeq _{B}$, that is, $\succcurlyeq $ is the partial order on $%
X:=A\sqcup B$ with $x\succcurlyeq y$ iff either $x\succcurlyeq _{A}y$ or $%
x\succcurlyeq _{B}y$. Now pick any number $\theta \geq \max \{$diam$(A),$diam%
$(B)\},$ and consider the function $D:X\times X\rightarrow \lbrack 0,\infty )
$ with%
\begin{equation*}
D(x,y):=\left\{ 
\begin{array}{ll}
d(x,y), & \text{if }(x,y)\in A^{2}\text{ or }(x,y)\in B^{2} \\ 
\frac{1}{2}\theta , & \text{otherwise.}%
\end{array}%
\right. 
\end{equation*}%
It is easily checked that $D$ is a metric on $X.$ In fact, $%
(X,D,\succcurlyeq )$ is a radial partially ordered metric space. $\square $

\bigskip 

\noindent \textit{Example 4. }Let $I$ stand for the unit interval $[0,1],$
and take any set $J$ that does not intersect $I.$ Define the partial order
on $X:=I\sqcup J$ with $x\succcurlyeq y$ iff either $(x,y)\in J\times I$ or $%
\{x,y\}\subseteq I$ and $x\geq y.$ (In other words, $\succcurlyeq $ agrees
with the usual order on $I,$ and puts anything in $J$ above all numbers in $%
I.$ No two distinct elements of $J$ are $\succcurlyeq $-comparable.) Define $%
d:X\times X\rightarrow \lbrack 0,\infty )$ as follows: (i) $d|_{I\times I}$
is the absolute value metric on $I$; (ii) $d|_{J\times J}$ is the discrete
metric on $J$; (iii) $d(x,y):=1+y$ if $(x,y)\in J\times I;$ and (iv) $%
d(x,y):=1+x$ if $(x,y)\in I\times J.$ Then, $(X,d,\succcurlyeq )$ is a
radial partially ordered metric space. $\square $

\bigskip

In passing, we note that it may be a mistake to think of the radiality
property as prohibitively strong. In the context of metric data analysis and
machine learning (see \cite{C-F-S}, \cite{D-Y-W-X}, and \cite{F-S}), one
often works with finite metric spaces or metric graphs (relative to which
the Lipschitz extension problems are by no means trivial). As witnessed by
Examples 1 and 2 above, the radiality property may turn out to be
considerably less demanding in those sorts of environments.

\subsection{Lipschitz Functions}

For any real number $K>0,$ a function $f:X\rightarrow Y$ from a partially
ordered metric space $X=(X,d_{X},\succcurlyeq _{X})$ to a partially ordered
metric space $Y=(Y,d_{Y},\succcurlyeq _{Y})$ is said to be $K$\textit{%
-Lipschitz} if for every $x,y\in X$, 
\begin{equation}
d_{Y}(f(x),f(y))\leq Kd_{X}(x,y).  \label{lip}
\end{equation}%
We say that $f$ is \textit{Lipschitz }if it is $K$-Lipschitz for some $K\geq
0.$ The smallest $K\geq 0$ such that (\ref{lip}) holds for every $x,y\in X$,
is called the \textit{Lipschitz constant} of $f.$ For excellent treatments
of the general theory of Lipschitz functions, see \cite{C-M-N} and \cite%
{Weaver}.

We denote the set of all $K$-Lipschitz maps from $X$ to $Y$ as Lip$_{K}(X,Y),
$ but write Lip$_{K}(X)$ for Lip$_{K}(X,\mathbb{R}).$ In turn, the sets of
all order-preserving members of Lip$_{K}(X,Y)$ and Lip$_{K}(X)$ are denoted
as Lip$_{K,\uparrow }(X,Y)$ and Lip$_{K,\uparrow }(X),$ respectively.
Throughout this note, we consider these as metric spaces relative to the
uniform metric. This makes these spaces complete, but in general not
separable.

\subsection{The Monotone Lipschitz Extension Property}

We say that a partially ordered metric space $(X,d,\succcurlyeq )$ has the 
\textit{monotone Lipschitz extension property }if for every nonempty $%
S\subseteq X$, every $K>0$ and $f\in $ Lip$_{K,\uparrow }(S),$ there exists
an $F\in $ Lip$_{K,\uparrow }(X)$ with $F|_{S}=f.$ In this terminology, the
classical \textit{McShane-Whitney extension theorem} can be viewed as saying
that $(X,d,=)$ has the monotone Lipschitz extension property. Our primary
objective in this note is to see exactly to what extent we can replace $=$
with a partial order on $X$ in this statement.

\bigskip

\noindent \textit{Remark 1. }When $(X,d,\succcurlyeq )$ has the monotone
Lipschitz extension property, we can always ensure the achieved extension
have the same range as the function to be extended. To see this, take any $%
F\in $ Lip$_{K,\uparrow }(X)$ and $S\subseteq X.$ Where $m:=\inf_{x\in S}F(x)
$ and $M:=\sup_{x\in S}F(x)$, the map $G:X\rightarrow \lbrack m,M]$ defined
by%
\begin{equation*}
G(x):=\max \{\min \{F(x),M\},m\},
\end{equation*}%
is an $\succcurlyeq $-increasing $K$-Lipschitz map with\textsl{\ }$%
G|_{S}=F|_{S}$. $\square $

\section{Monotone Lipschitz Extensions}

Unless a partially ordered metric space is totally ordered, or it is finite,
its radiality seems like a fairly demanding condition. However, our main
finding in this note shows that this condition is necessary and sufficient
for any such space to possess the monotone Lipschitz extension property.

\bigskip

\noindent \textsf{\textbf{Theorem 2.}}\textbf{\ }\textit{A partially ordered
metric space }$(X,d,\succcurlyeq )$\textit{\ has the monotone Lipschitz
extension property if and only if it is radial.}

\medskip

\textit{\textbf{Proof.}}\textsf{\textbf{\ }}Suppose $(X,d,\succcurlyeq )$ is
not radial. Then, there exist three points $x,y,z$ in $X$ such that either 
\begin{equation}
x\succcurlyeq ^{\bullet }y\succ z\hspace{0.2in}\text{and\hspace{0.2in}}%
d(x,z)<d(x,y),  \label{w1}
\end{equation}%
or%
\begin{equation}
x\succ y\succcurlyeq ^{\bullet }z\hspace{0.2in}\text{and\hspace{0.2in}}%
d(x,z)<d(y,z)\text{.}  \label{w2}
\end{equation}%
Assume first the case (\ref{w1}), set $S:=\{x,y\},$ and define $%
f:S\rightarrow \mathbb{R}$ by $f(x):=d(x,y)$ and $f(y):=0.$ Then, $f\in $ Lip%
$_{1,\uparrow }(S),$ but for any 1-Lipschitz extension $F:X\rightarrow 
\mathbb{R}$ of $f,$ we have%
\begin{equation*}
F(z)\geq F(x)-d(x,z)>f(x)-d(x,y)=0=F(y)
\end{equation*}%
which means $F$ is not $\succcurlyeq $-increasing. If, on the other hand, (%
\ref{w2}) holds, we set $S:=\{y,z\},$ and define $f:S\rightarrow \mathbb{R}$
by $f(y):=d(y,z)$ and $f(z):=0.$ Then, $f\in $ Lip$_{1,\uparrow }(S),$ but
for any 1-Lipschitz extension $F:X\rightarrow \mathbb{R}$ of $f,$ we have%
\begin{equation*}
F(x)\leq F(z)+d(x,z)<d(y,z)=F(y)
\end{equation*}%
which means $F$ is not $\succcurlyeq $-increasing. This proves the necessity
part of the assertion. The sufficiency part is a special case of a more
general result that we will establish below. $\blacksquare $

\bigskip

There does not seem like there is an easy way of getting around the
radiality requirement for the monotonic Lipschitz extension problem. For a
partially ordered metric space $(X,d,\succcurlyeq )$ that is not radial, the
argument above shows that it may not be possible to extend an $\succcurlyeq $%
-increasing 1-Lipschitz map on a compact and $\succcurlyeq $-increasing (or $%
\succcurlyeq $-decreasing) set $S\subseteq X$ to an $\succcurlyeq $%
-increasing 1-Lipschitz map on $X.$

Setting $\succcurlyeq $ as the equality relation in Theorem 2 yields the
classical McShane-Whitney extension theorem. The following is another
straightforward corollary.

\bigskip

\noindent \textsf{\textbf{Corollary 3.}}\textbf{\ }\textit{A linearly
ordered metric space }$(X,d,\succcurlyeq )$\textit{\ has the monotone
Lipschitz extension property if and only if it is radially convex.}

\bigskip

\noindent \textit{Remark 2. }It was shown by Mehta \cite{Min} that every
topological loset $(X,\succcurlyeq )$ is a normally ordered topological
space. Therefore, specializing the Nachbin extension theorem to the context
of metric spaces, we find: \textit{Given any metric loset }$%
(X,d,\succcurlyeq )$\textit{, and any compact }$S\subseteq X,$\textit{\
every }$\succcurlyeq $\textit{-increasing }$f\in C(S)$\textit{\ extends to
an }$\succcurlyeq $\textit{-increasing }$F\in C(X).$ Corollary 3 can be
thought of as the reflection of this result in the context of Lipschitz
functions. As a return to adding the hypothesis of radial convexity to the
picture, it achieves an order-preserving Lipschitz extension of any
order-preserving Lipschitz function defined on any (possibly non-compact)
subset of $X.$ $\square $

\bigskip

The Lipschitz extension problem for Banach space-valued maps on a metric
space is a rather deep one, and is the subject of ongoing research in metric
space theory and geometric functional analysis. However, there is one
special case of the problem which is settled by the McShane-Whitney theorem
in a routine manner. This is when the Lipschitz maps to be extended take
values in the Banach space $\ell _{\infty }(T)$ of all bounded real
functions on some nonempty set $T.$ (This generalization is of interest,
because every metric space can be isometrically embedded in $\ell _{\infty
}(T)$ for some $T$.) Precisely the same holds for the monotone Lipschitz
extension problem as well where we consider $\ell _{\infty }(T)$ as
partially ordered coordinatewise. (For any $u,v\in \ell _{\infty }(T)$, we
write $u\geq v$ whenever $u(t)\geq v(t)$ for every $t\in T$). We now prove
the sufficiency part of Theorem 2 in this more general context.

\bigskip

\noindent \textsf{\textbf{Theorem 4.}}\textbf{\ }\textit{Let }$%
(X,d,\succcurlyeq )$\textit{\ be a radial partially ordered metric space.
For any }$K\geq 0,$\textit{\ let }$S$\textit{\ be a nonempty subset of }$X$%
\textit{\ and }$f:S\rightarrow \ell _{\infty }(T)$\textit{\ an
order-preserving }$K$\textit{-Lipschitz map. Then, there exists an
order-preserving }$K$\textit{-Lipschitz map }$F:X\rightarrow \ell _{\infty
}(T)$\textit{\ with }$F|_{S}=f$\textit{.}

\medskip

\textit{\textbf{Proof.}}\textsf{\textbf{\ }}We assume $S\neq X,$ for
otherwise there is nothing to prove. Similarly, the claim is trivially true
when $K=0,$ so we may assume $K>0.$ Moreover, it is enough to prove the
assertion for $K=1$, for then the general case obtains by applying what is
established to the map $\frac{1}{K}f.$

The following proof is patented after the typical way one proves the
Hahn-Banach Theorem. In the initial stage of the argument, we take an
arbitrary $x\in X\backslash S$ and extend $f$ to an order-preserving
1-Lipschitz function on $S\cup \{x\}$. To this end, consider the functions $%
a_{x}:T\rightarrow \lbrack -\infty ,\infty ]$ and $b_{x}:T\rightarrow
\lbrack -\infty ,\infty ]$ defined as 
\begin{equation*}
a_{x}(t):=\sup \left\{ f(z)(t):z\in S\cap x^{\downarrow }\right\} 
\end{equation*}%
and%
\begin{equation*}
b_{x}(t):=\inf \left\{ f(y)(t):y\in S\cap x^{\uparrow }\right\} \text{.}
\end{equation*}%
If $S\cap x^{\downarrow }=\varnothing ,$ then $a_{x}(t)=-\infty $ for every $%
t\in T,$ while $S\cap x^{\uparrow }=\varnothing $ implies $b_{x}(t)=\infty $
for every $t\in T$. On the other hand, if both $S\cap x^{\downarrow }$ and $%
S\cap x^{\uparrow }$ are nonempty, monotonicity of $f$ yields $-\infty
<a_{x}(t)\leq b_{x}(t)<\infty $ for all $t\in T.$ In all contingencies,
then, $[a_{x}(t),b_{x}(t)]$ is a nonempty interval in the set of all
extended reals.

We next define the functions $\alpha _{x}:T\rightarrow \lbrack -\infty
,\infty ]$ and $\beta _{x}:T\rightarrow \lbrack -\infty ,\infty ]$ by 
\begin{equation*}
\alpha _{x}(t):=\sup \left\{ f(z)(t)-d(x,z):z\in S\right\} 
\end{equation*}%
and%
\begin{equation*}
\beta _{x}(t):=\inf \left\{ f(y)(t)+d(x,y):y\in S\right\} \text{.}
\end{equation*}%
(These are the McShane and Whitney extensions of $f,$ respectively.) In this
case, both $\alpha _{x}(t)$ and $\beta _{x}(t)$ are real numbers for every $%
t\in T.$ In fact, as $f$ is 1-Lipschitz, for every $y,z\in S$ we have 
\begin{equation*}
f(z)(t)-f(y)(t)\leq \left\Vert f(z)-f(y)\right\Vert _{\infty }\leq
d(z,y)\leq d(x,y)+d(x,z),
\end{equation*}%
whence $f(z)(t)-d(x,z)\leq f(y)(t)+d(x,y)$, for all $t\in T.$ Conclusion: $%
-\infty <\alpha _{x}(t)\leq \beta _{x}(t)<\infty $ for all $t\in T.$

We claim that 
\begin{equation}
\alpha _{x}(t)\leq b_{x}(t)\hspace{0.2in}\text{and\hspace{0.2in}}%
a_{x}(t)\leq \beta _{x}(t)  \label{main}
\end{equation}%
for every $t\in T.$ To see this, suppose $\alpha _{x}(t)>b_{x}(t)$ for some $%
t\in T.$ Then, there exist $y\in S\cap x^{\uparrow }$ and $z\in S$ such that 
$f(y)(t)<f(z)(t)-d(x,z).$ It follows that $f(y)(t)<f(z)(t),$ so $%
y\succcurlyeq z$ does not hold (because $f$ is $\succcurlyeq $-increasing).
Thus: $z\succcurlyeq ^{\bullet }y\succ x.$ Since $(X,d,\succcurlyeq )$ is
radial, therefore, $d(x,z)\geq d(y,z).$ This entails%
\begin{equation*}
f(y)(t)<f(z)(t)-d(x,z)\leq f(z)(t)-d(y,z),
\end{equation*}%
and hence, $\left\Vert f(z)-f(y)\right\Vert _{\infty }\geq
f(z)(t)-f(y)(t)>d(z,y),$ contradicting $f$ being 1-Lipschitz. We conclude
that $\alpha _{x}(t)\leq b_{x}(t)$ for all $t\in T,$ as claimed. The second
inequality in (\ref{main}) is established analogously.

In view of these observations, we conclude that the intervals $%
[a_{x}(t),b_{x}(t)]$ and $[\alpha _{x}(t),\beta _{x}(t)]$ overlap for every $%
t\in T.$ We define $F:S\cup \{x\}\rightarrow \ell _{\infty }(T)$ as%
\begin{equation*}
F(w)(t):=\left\{ 
\begin{array}{ll}
f(w)(t), & \text{if }w\in S \\ 
\theta (t), & \text{if }w=x,%
\end{array}%
\right. 
\end{equation*}%
where $\theta (t)$ is an arbitrarily picked real number in $%
[a_{x}(t),b_{x}(t)]\cap \lbrack \alpha _{x}(t),\beta _{x}(t)]$ for any $t\in
T.$ Then, $F$ is 1-Lipschitz, because for any $y\in S$, we have%
\begin{equation*}
f(y)(t)-d(x,y)\leq \alpha _{x}(t)\leq F(x)(t)\leq \beta _{x}(t)\leq
f(y)(t)+d(x,y),
\end{equation*}%
and hence $\left\vert F(x)(t)-F(y)(t)\right\vert \leq d(x,y),$ for all $t\in
T,$ that is, $\left\Vert F(y)-F(x)\right\Vert _{\infty }$ $\leq d(x,y).$ On
the other hand, for every $y\in S$ with $y\succcurlyeq x,$ we have $%
f(y)(t)\geq b_{x}(t)\geq F(x)(t),$ and similarly, for every $z\in S$ with $%
x\succcurlyeq z,$ we have $f(x)(t)\geq a_{x}(t)\geq F(z)(t),$ for all $t\in
T.$ Thus, $F$ is order-preserving as well.

The proof is completed by a standard transfinite induction argument. Let $%
\mathcal{F}$ stand for the set of all $(A,F)$ such that $S\subseteq
A\subseteq X$ and $F\in $ Lip$_{1,\uparrow }(A,\ell _{\infty }(T))$ with $%
F|_{S}=f.$ Since it includes $(S,f),$ this collection is not empty. In
addition, it is easily verified that $(\mathcal{F},\trianglerighteq )$ is an
inductive poset where $(A,F)\trianglerighteq (B,G)$ iff $A\supseteq B$ and $%
F|_{B}=G.$ So, by Zorn's Lemma, there is a $\trianglerighteq $-maximal
element $(A,F)$ in $\mathcal{F}$. In view of the first part of the proof, we
must have $T=X.$ $\blacksquare $

\bigskip

\noindent \textit{Remark 3. }By setting $\theta (t):=\max \{a_{x}(t),\alpha
_{x}(t)\}$ for all $t\in T$ in the proof above, and modifying the
transfinite induction part of the proof in the obvious way, we find that
there is a smallest order-preserving $K$-Lipschitz map $F:X\rightarrow \ell
_{\infty }(T)$ with $F|_{S}=f$ in the context of Theorem 4. That there is
also a largest such $F$ is established analogously. $\square $

\bigskip

\noindent \textit{Remark 4. }There are various generalizations of the
Lipschitz property, and the construction above adapts to some of these. To
wit, Miculescu \cite{Mic} considers $(K,g)$-Lipschitz functions which are
functions $f$ from a metric space $(X,d_{X})$ to another metric space $%
(Y,d_{Y})$ such that $d_{X}(f(x),f(y))\leq Kd_{Y}(g(x),g(y))$ for every $%
x,y\in X.$ Theorems 2 and 4 modify in the obvious way to account for such
functions as well. $\square $

\bigskip

\noindent \textit{Remark 5. }Given Theorem 2, it is natural to inquire if
the monotonic Lipschitz extensions of real functions can be carried out
locally. To state the problem, we recall that a real map on a metric space $%
X=(X,d)$ is called \textit{pointwise Lipschitz} if for every $y\in X,$ there
exist $K_{y}\geq 0$ and $\delta _{y}>0$ such that $\left\vert
f(x)-f(y)\right\vert \leq K_{y}d(x,y)$ for all $x\in X$ with $d(x,y)<\delta
_{y}.$ The question is: If $(X,d,\succcurlyeq )$ is a radial partially
ordered metric space, $S$ a nonempty closed subset of $X$ and $%
f:S\rightarrow \mathbb{R}$ is an $\succcurlyeq $-increasing pointwise
Lipschitz map, does there exist an $\succcurlyeq $-increasing pointwise
Lipschitz map $F:X\rightarrow \mathbb{R}$ with $F|_{S}=f$? If $\succcurlyeq $
is the equality relation, the answer is known to be yes; see, for instance, 
\cite{C-G} and \cite{Gutev}. The first part of the proof above also adapts
to show that the answer is yes so long as we add only finitely many points
in the extension. That is, minor modifications of that part of the proof
yields the following fact:

\medskip

\noindent \textit{Let }$(X,d,\succcurlyeq )$\textit{\ be a radial partially
ordered metric space, and }$S$\textit{\ a nonempty closed subset of }$X$%
\textit{\ with }$\left\vert X\backslash S\right\vert <\infty .$\textit{\
Then, every }$\succcurlyeq $\textit{-increasing pointwise Lipschitz map on }$%
S$\textit{\ can be extended to an }$\succcurlyeq $\textit{-increasing
pointwise Lipschitz map on }$X$\textit{.}

\medskip

\noindent Unfortunately, the transfinite inductive step of the proof above
fails to deliver this result without the requirement $\left\vert X\backslash
S\right\vert <\infty .$ $\square $

\section{Functional Representations of Radial Orders}

For any nonempty $X$, $\mathcal{F}\subseteq \mathbb{R}^{X},$ and $x,y\in X,$
we write $\mathcal{F}(x)\geq \mathcal{F}(y)$ to mean $f(x)\geq f(y)$ for
every $f\in \mathcal{F}$. For any such collection $\mathcal{F},$ the binary
relation $\succsim $ on $X$ defined by $x\succsim y$ iff $\mathcal{F}(x)\geq 
\mathcal{F}(y),$ is a preorder on $X.$ Conversely, for every preorder $%
\succsim $ on $X,$ there is a family $\mathcal{F}$ with $x\succsim y$ iff $%
\mathcal{F}(x)\geq \mathcal{F}(y)$ for every $x,y\in X.\footnote{%
This is readily proved by taking $\mathcal{F}$ as the set of all indicator
functions of $\{z\in X:z\succsim x\}$ as $x$ varies over $X$.}$ In this
case, we say that $\mathcal{F}$ \textit{represents} $\succsim $. In several
applied mathematical fields, such as decision theory and the theory of
optimal transportation, it is important to determine the structure of the
families of real functions that may represent a given preorder in this sense.

As an easy consequence of Theorem 2, we find that any radial partial order
on any metric space can be represented by a family of order-preserving
1-Lipschitz real-valued functions.

\bigskip

\noindent \textsf{\textbf{Proposition 5.}}\textbf{\ }\textit{Let }$%
(X,d,\succcurlyeq )$\textit{\ be a radial partially ordered metric space.
Then, there exists an }$\mathcal{F}\subseteq $\textit{\ }Lip$_{1,\uparrow
}(X)$ \textit{that represents} $\succcurlyeq $. \textit{If }$(X,d)$\textit{\
is compact, we can choose }$\mathcal{F}$\textit{\ in such a way that it is
compact and} $\sup_{F\in \mathcal{F}}\left\Vert F\right\Vert _{\infty }\leq $
diam$(X)$.

\medskip

\textit{\textbf{Proof.}}\textsf{\textbf{\ }}Assume $\left\vert X\right\vert
>1,$ which implies $\succcurlyeq ^{\bullet }$ $\neq \varnothing ,$ for
otherwise there is nothing to prove. For any $x,y\in X$ with $x\succcurlyeq
^{\bullet }y,$ define $f_{x,y}\in \mathbb{R}^{\{x,y\}}\mathbb{\ }$by $%
f_{x,y}(x):=d(x,y)$ and $f_{x,y}(y):=0,$ and note that $f\in $ Lip$%
_{1,\uparrow }(\{x,y\}).$ We apply Theorem 2 to extend $f_{x,y}$ to an $%
\succcurlyeq $-increasing 1-Lipschitz real-valued map $F_{x,y}$ on $X.$
Next, define $\mathcal{F}:=\{F_{x,y}:x\succcurlyeq ^{\bullet }y\}.$ Then, $%
x\succcurlyeq y$ implies $F(x)\geq F(y)$ for all $F\in \mathcal{F}$ simply
because every member of $\mathcal{F}$ is $\succcurlyeq $-increasing.
Conversely, if $x\succcurlyeq y$ does not hold, we have $F(x)<F(y)$ for some 
$F\in \mathcal{F}$, namely, $F=F_{y,x}$.

Now suppose $(X,d)$ is compact. Put $K:=$ diam$(X),$ and note that $K\in
(0,\infty ).$ Next, for any fixed $e\in X,$ define 
\begin{equation*}
\mathcal{G}:=\{\tfrac{1}{K}(F-F(e)):F\in \mathcal{F}\}.
\end{equation*}%
Then, $\left\vert G(x)\right\vert =\left\vert G(x)-G(e)\right\vert \leq
K^{-1}d(x,e)\leq 1$ for every $x\in X,$ so $\left\Vert G\right\Vert _{\infty
}\leq 1,$ for every $G\in \mathcal{G}$. Moreover, $\mathcal{G}\subseteq $ Lip%
$_{1/K,\uparrow }(X)$ and $\mathcal{G}$ represents $\succcurlyeq $. Then, $%
\mathcal{H}:=$ $K$cl$(\mathcal{G)}$ is a closed and bounded set of $1$%
-Lipschitz bounded functions which represents $\succcurlyeq $. Since any
subset of Lip$_{1}(X)$ is equicontinuous, applying the Arzel\`{a}-Ascoli
Theorem yields the second claim of the proposition. $\square $

\bigskip 

As an immediate consequence of Proposition 5 we obtain the somewhat
surprising fact that every radial partially ordered metric space is, per
force, a metric poset. That is:

\bigskip

\noindent \textsf{\textbf{Corollary 6.}}\textbf{\ }\textit{Every radial
partial order on a metric space }$X$\textit{\ is a closed subset of }$%
X\times X$\textit{.}

\bigskip 

The concepts of \textquotedblleft radial partially ordered metric
space\textquotedblright\ and \textquotedblleft radial metric
poset\textquotedblright\ are thus identical. We will adopt the latter
terminology in the remainder of the paper. 

We next apply our main extension theorem to show that a radially convex
linear order on a $\sigma $-compact metric space can be represented by a
Lipschitz function. The main ingredient of the argument is contained in the
following observation.

\bigskip

\noindent \textsf{\textbf{Lemma 7.}}\textbf{\ }\textit{Let }$%
(X,d,\succcurlyeq )$\textit{\ be a radial metric poset. Then, for any
compact subset }$S$\textit{\ of }$X,$ \textit{there is an }$G\in $ Lip$%
_{1,\uparrow }(X)$\textit{\ such that }$\left\Vert G\right\Vert _{\infty
}\leq $ diam$(S)$ and%
\begin{equation*}
G(x)>G(y)\hspace{0.2in}\text{\textit{for every }}x,y\in S\text{ \textit{with}
}x\succ y.
\end{equation*}

\textit{\textbf{Proof.}}\textsf{\textbf{\ }}Take any compact $S\subseteq X,$
and use Proposition 5 to find a compact, and hence separable, $\mathcal{F}%
\subseteq $ Lip$_{1,\uparrow }(S)$ such that (i) $\sup_{F\in \mathcal{F}%
}\left\Vert F\right\Vert _{\infty }\leq $ diam$(S);$ and (ii) $x\succcurlyeq
y$ iff $\mathcal{F}(x)\geq \mathcal{F}(y)$ for every $x,y\in X.$ Let $(F_{m})
$ be a sequence in $\mathcal{F}$ such that $\{F_{1},F_{2},...\}$ is dense in 
\QTR{cal}{F}. We define $G:=\sum_{n\geq 1}2^{-n}F_{n}$. It is readily
checked that $F\in $ Lip$_{1,\uparrow }(S).$ Besides, if $x,y\in S$ satisfy $%
x\succ y,$ then $F(x)>F(y)$ for some $F\in \mathcal{F}$ (because $\mathcal{F}
$ represents $\succcurlyeq $). Consequently, since $\{F_{1},F_{2},...\}$ is
dense in $\mathcal{F}$ relative to the uniform metric, there exists an $n\in 
\mathbb{N}$ with $F_{n}(x)>F_{n}(y),$ which implies $G(x)>G(y).$ To complete
the proof, we extend $G$ to $X$ by using Theorem 2, and recall Remark 1. $%
\square $

\bigskip

\noindent \textsf{\textbf{Theorem 8.}}\textbf{\ }\textit{Let }$%
(X,d,\succcurlyeq )$\textit{\ be a radial metric poset such that }$(X,d)$%
\textit{\ is }$\sigma $\textit{-compact. Then, there is a Lipschitz function 
}$F:X\rightarrow \mathbb{R}$\textit{\ with }%
\begin{equation}
F(x)>F(y)\hspace{0.2in}\text{\textit{for every }}x,y\in S\text{ \textit{with}
}x\succ y.  \label{HH}
\end{equation}

\textit{\textbf{Proof.}}\textsf{\textbf{\ }}By hypothesis, there exists a
sequence $(S_{m})$ of compact subsets of $X$ such that $S_{1}\subseteq
S_{2}\subseteq \cdot \cdot \cdot $ and $S_{1}\cup S_{2}\cup \cdot \cdot
\cdot =X.$ We may assume $\left\vert S_{1}\right\vert >1.$ Put $K_{n}:=$ diam%
$(S_{n}),$ and note that $K_{n}\in (0,\infty )$ for each $n.$ By Lemma 7,
for every $n\in \mathbb{N},$ there is a $G_{n}\in $ Lip$_{1,\uparrow }(X)$%
\textit{\ }such that\textit{\ }$\left\Vert G_{n}\right\Vert _{\infty }\leq
K_{n}$ and $G_{n}(x)>G_{n}(y)$ for every $x,y\in S_{n}$ with $x\succ y.$ We
define $F\in \mathbb{R}^{X}$ by $F(x):=\sum_{n\geq 1}2^{-n}K_{n}^{-1}G_{n}.$
It is plain that $F(x)>F(y)$ for every $x,y\in X$ with $x\succ y.$ Moreover, 
$F$ is $K_{1}^{-1}$-Lipschitz. Indeed for any $x,y\in X,$%
\begin{equation*}
\left\vert F(x)-F(y)\right\vert \leq \sum_{n\geq 1}\tfrac{1}{2^{n}K_{n}}%
\left\vert G_{n}(x)-G_{n}(y)\right\vert \leq \sum_{n\geq 1}\tfrac{1}{%
2^{n}K_{1}}d(x,y)\leq \tfrac{1}{K_{1}}d(x,y)
\end{equation*}%
since $K_{1}\leq K_{n}$ for each $n$. $\square $

\bigskip

\noindent \textsf{\textbf{Corollary 9.}}\textbf{\ }\textit{Let }$%
(X,d,\succcurlyeq )$\textit{\ be a radially convex metric loset such that }$%
(X,d)$\textit{\ is }$\sigma $-\textit{compact. Then, there is a Lipschitz
function }$F:X\rightarrow \mathbb{R}$ \textit{with }%
\begin{equation*}
x\succcurlyeq y\hspace{0.2in}\text{\textit{if and only if\hspace{0.2in}}}%
F(x)\geq F(y)
\end{equation*}%
\textit{for every }$x,y\in X.$

\bigskip

This result has the flavor of continuous utility representation theorems of
decision theory. Indeed, it provides a rather easy proof of the following
well-known result of that literature.

\bigskip

\noindent \textsf{\textbf{Corollary 10.}}\textbf{\ }\textit{Let }$\succsim $ 
\textit{be a closed total preorder on a compact metric space }$X=(X,d)$%
\textit{. Then, there exists a continuous map }$u:X\rightarrow \mathbb{R}$%
\textit{\ such that }%
\begin{equation*}
x\succsim y\hspace{0.2in}\text{\textit{if and only if\hspace{0.2in}}}%
u(x)\geq u(y)
\end{equation*}%
\textit{for every }$x,y\in X.$

\medskip

\textit{\textbf{Proof.}}\textsf{\textbf{\ }}Define $\mathbf{x}:=\{y\in
X:x\succsim y\succsim x\}$ for any $x\in X,$ and note that $\mathcal{X}:=\{%
\mathbf{x}:x\in X\}$ is a partition of $X.$ Then, the binary relation $%
\succcurlyeq $ $\subseteq \mathcal{X}\times \mathcal{X}$ defined by $\mathbf{%
x}\succcurlyeq \mathbf{y}$ iff $x\succsim y,$ is a partial order on $%
\mathcal{X}.$ Let $H_{d}$ stand for the Hausdorff metric on $\mathcal{X}$.
Then, $(\mathcal{X},H_{d},\succcurlyeq )$ is a compact metric loset. By the
Carruth metrization theorem (of \cite{Carruth}), there exists a metric $D$
on $\mathcal{X}$ such that $H_{d}$ and $D$ are equivalent, and $D(\mathbf{x,z%
})=D(\mathbf{x,y})+D(\mathbf{y,z})$ for every $x,y,z\in X$ with $x\succ
y\succ z.$ We may thus apply Corollary 9 to obtain an $\succcurlyeq $%
-increasing and 1-Lipschitz $F$ map on $(\mathcal{X},D,\succcurlyeq )$ such
that $\mathbf{x}\succcurlyeq \mathbf{y}$ iff $F(\mathbf{x})\geq F(\mathbf{y})
$ for every $x,y\in X.$ The map $u:X\rightarrow \mathbb{R}$ with $u(x):=F(%
\mathbf{x})$ fulfills the requirements of the assertion. $\square $

\section{Monotone Uniformly Continuous Extensions}

\subsection{A Monotone Version of McShane's Uniformly Continuous Extension
Theorem}

As an another application of Theorem 2, we prove a uniformly continuous
extension theorem in the context of radial metric posets. A special case of
this theorem will correspond to the monotonic version of McShane's famous
uniformly continuous extension theorem for bounded functions.

For any metric spaces $X=(X,d_{X})$ and $Y=(Y,d_{Y}),$ a function $%
f:X\rightarrow Y$ is said to be \textit{Lipschitz for large distances }if
for every $\delta >0$ there is a $K_{\delta }>0$ such that $%
d_{Y}(f(x),f(y))\leq K_{\delta }d_{X}(x,y)$ whenever $d_{X}(x,y)\geq \delta .
$ This concept often arises with extension and approximation problems
concerning uniformly continuous functions; see, for instance, \cite{L-R}, 
\cite{G-J} and \cite{B-L}. In fact, a basic result of this literature says
that every uniformly continuous map on a Menger-convex metric space is, per
force, Lipschitz for large distances (cf. \cite[Proposition 1.11]{B-L}).

We need to make two observations about real-valued functions that are
Lipschitz for large distances. The first one is basic, and was noted
explicitly in \cite{G-J}.

\bigskip 

\noindent \textsf{\textbf{Lemma 11.}}\textbf{\ }\textit{Every bounded
real-valued function on a metric space is Lipschitz for large distances.}

\medskip

\textit{\textbf{Proof.}}\textsf{\textbf{\ }}For any bounded real-valued
function $f$ on a metric space $X=(X,d),$ and $\delta >0,$ we have%
\begin{equation*}
\left\vert f(x)-f(y)\right\vert \leq \left( \tfrac{2\left\Vert f\right\Vert
_{\infty }}{\delta }\right) d(x,y)
\end{equation*}%
for all $x,y\in X$ with $d(x,y)\geq \delta .$ $\blacksquare $

\bigskip 

Our second observation provides a characterization of uniformly continuous
real-valued maps that are Lipschitz for large distances. This
characterization seems new, but we should note that Beer and Rice \cite{B-R}
work out several related results. In the statement of the result, and
henceforth, $\omega _{f}$ stands for the \textit{modulus of continuity }of
any given real-valued function $f$ on $X=(X,d),$ that is, $\omega
_{f}:[0,\infty )\rightarrow \lbrack 0,\infty ]$ is the function defined by%
\begin{equation*}
\omega _{f}(t):=\sup \{\left\vert f(x)-f(y)\right\vert :x,y\in X\text{ and }%
d(x,y)\leq \delta \}.
\end{equation*}

\medskip 

\noindent \textsf{\textbf{Lemma 12.}}\textbf{\ }\textit{Let }$X=(X,d)$ 
\textit{be a metric space and }$f\in UC(X).$ \textit{Then, }$f$ \textit{is
Lipschitz for large distances if and only if there exist nonnegative real
numbers }$a$\textit{\ and }$b$\textit{\ such that }$\omega _{f}(t)\leq at+b$%
\textit{\ for every }$t\geq 0.$\footnote{\cite{B-R} refers to an $%
f:X\rightarrow \mathbb{R}$ with the latter property as a function \textit{%
having an affine majorant}, and investigates it in detail. In fact, this
concept already plays a prominent role in McShane's original article \cite%
{McS} where (on its page 841) it is emphasized that when $X$ is a normed
linear space, a uniformly extendable real function on a subset of $X$ must
have an affine majorant. }

\medskip

\textit{\textbf{Proof.}}\textsf{\textbf{\ }}For any $a,b\in \mathbb{R},$ let 
$h_{a,b}$ denote the map $t\mapsto at+b$ on $[0,\infty )$. Suppose first
that $\omega _{f}\leq h_{a,b}$ for some $a,b\geq 0.$ Then, for any $\delta
>0,$ setting $K_{\delta }:=a+b/\delta $ yields 
\begin{equation*}
\left\vert f(x)-f(y)\right\vert \leq \omega _{f}(d(x,y))\leq ad(x,y)+b\leq
ad(x,y)+b\left( \tfrac{d(x,y)}{\delta }\right) =K_{\delta }d(x,y)
\end{equation*}%
for every $x,y\in X$ with $d(x,y)\geq \delta .$ Conversely, suppose $f$ is
Lipschitz for large distances. Note first that uniform continuity of $f$
entails that there is a $\delta >0$ with $\omega _{f}(\delta )\leq 1.$ In
turn, by the Lipschitz property of $f,$ there exists a $K:=K_{\delta }>0$
such that 
\begin{equation*}
\left\vert f(x)-f(y)\right\vert \leq Kd(x,y)\text{\hspace{0.2in}for all }%
x,y\in X\text{ with }d(x,y)\geq \delta .
\end{equation*}
We wish to show that $\omega _{f}\leq h_{K,1}.$ To this end, fix an
arbitrary $t\geq 0,$ and take any $x,y\in X$ with $d(x,y)\leq t.$ If $%
d(x,y)<\delta ,$ then $\left\vert f(x)-f(y)\right\vert \leq \omega
_{f}(\delta )\leq 1\leq h_{K,1}(t).$ Otherwise, $\left\vert
f(x)-f(y)\right\vert \leq Kd(x,y)\leq Kt\leq h_{K,1}(t).$ Conclusion: $%
\left\vert f(x)-f(y)\right\vert \leq h_{K,1}(t)$ for any $x,y\in X$ with $%
d(x,y)\leq t.$ Taking the sup over all such $x$ and $y$ yields $\omega
_{f}(t)\leq h_{K,1}(t).$ $\blacksquare $

\bigskip 

We now proceed to show that the uniformly continuous extension theorem of
McShane \cite{McS} also generalizes to the context of radial metric posets.
This is proved most easily by adopting the remetrization technique of Beer 
\cite[pp. 23-25]{BeerBook} which derives the said extension from the
McShane-Whitney theorem. For the sake of completeness, we provide the
details of Beer's technique within the proof.

\bigskip

\noindent \textsf{\textbf{Theorem 13.}}\textbf{\ }\textit{Let }$%
(X,d,\succcurlyeq )$\textit{\ be a radial metric poset and }$S$ \textit{a
subset of }$X.$\textit{\ Then, for every }$\succcurlyeq $\textit{-increasing 
}$f\in UC(S)$ \textit{which is Lipschitz for large distances, there exists
an }$\succcurlyeq $\textit{-increasing }$F\in UC(X)$\textit{\ with }$F|_{S}=f
$\textit{.}

\medskip

\textit{\textbf{Proof.}}\textsf{\textbf{\ }}Let $\mathcal{H}$ stand for the
set of all increasing affine self-maps $h$ on $[0,\infty )$ with $\omega
_{f}\leq h.$ By Lemma 12, $\mathcal{H}\neq \varnothing $. We may thus define
the map $\varphi :[0,\infty )\rightarrow \mathbb{R}$ by $\varphi
(t):=\inf_{h\in \mathcal{H}}h(t).$ Clearly, $\varphi $ is an increasing and
concave (hence subadditive) self-map on $[0,\infty ).$ Since $f$ is not
constant, $\varphi (t)>0$ for some $t>0,$ so concavity of $\varphi $ entails 
$\varphi (t)>0$ for all $t>0.$ We claim that $\varphi $ is continuous at 0
(whence $\varphi \in C([0,\infty ))$) and $\varphi (0)=0$. To prove this,
take any $\varepsilon >0.$ Since $f$ is uniformly continuous, there exists a 
$\delta >0$ with $\omega _{f}(t)\leq \varepsilon $ for every $t\in \lbrack
0,\delta ).$ In turn, as $f$ is Lipschitz for large distances, there exists
a $K>0$ such that $\left\vert f(x)-f(y)\right\vert \leq Kd(x,y)$ whenever $%
d(x,y)\geq \delta .$ Now consider the self-map $h$ on $[0,\infty )$ with $%
h(t):=Kt+\varepsilon .$ Clearly, $\omega _{f}(t)\leq h(0)\leq h(t)$ for all $%
t\in \lbrack 0,\delta ),$ while $\omega _{f}(t)\leq \max \{\varepsilon
,Kt\}\leq h(t)$ for all $t\geq \delta .$ It follows that $h\in \mathcal{H}$.
But then $\varphi (t)\leq Kt+\varepsilon $ for all $t\geq 0,$ which implies $%
\inf_{t>0}\varphi (t)\leq \varepsilon .$ In view of the arbitrary choice of $%
\varepsilon ,$ we conclude that $\inf_{t>0}\varphi (t)=0=\varphi (0).$

With these preparations in place, we now turn to the task at hand. Define $%
D:X\times X\rightarrow \mathbb{R}$ by $D(x,y):=\varphi (d(x,y)).$ Since $%
\varphi (t)>0$ for all $t>0,$ it is obvious that $D(x,y)>0$ for every
distinct $x,y\in X,$ while $\varphi (0)=0$ implies $D(x,x)=0$ for all $x\in X
$. Moreover, $D$ is clearly symmetric and it satisfies the triangle
inequality (because $\varphi $ is increasing and subadditive). Thus: $%
(X,D,\succcurlyeq )$ is a partially ordered metric space. As $\varphi $ is
increasing, this space is radial. Besides, $\left\vert f(x)-f(y)\right\vert
\leq h(d(x,y))$ for every $x,y\in S$ and $h\in \mathcal{H},$ and it follows
that $\left\vert f(x)-f(y)\right\vert \leq D(x,y)$ for every $x,y\in S,$
that is, $f$ is 1-Lipschitz on the metric space $(S,D|_{S\times S}).$ By
Theorem 2, therefore, there exists an $\succcurlyeq $-increasing $%
F:X\rightarrow \mathbb{R}$ which is 1-Lipschitz on $(X,D)$ with $F|_{S}=f.$
But then for every $\varepsilon >0,$ continuity of $\varphi $ at $0$ ensures
that there is a $\delta >0$ small enough that $\varphi (t)<\varepsilon $ for
all $t\in (0,\delta ),$ which means $\left\vert F(x)-F(y)\right\vert
<\varepsilon $ for all $x,y\in X$ with $d(x,y)\leq \delta .$ It follows that 
$F$ is uniformly continuous on the metric space $(X,d).$ $\blacksquare $

\bigskip 

Since every bounded map on a metric space is Lipschitz for large distances
(Lemma 11), the following is a special case of Theorem 13. When $%
\succcurlyeq $ is taken as the equality relation in its statement, this
result reduces to McShane's uniformly continuous extension theorem for
bounded real-valued functions.

\bigskip 

\noindent \textsf{\textbf{Corollary 14.}}\textbf{\ }\textit{Let }$%
(X,d,\succcurlyeq )$\textit{\ be a radial metric poset. Every }$\succcurlyeq 
$\textit{-increasing, bounded and uniformly continuous map on a subset }$S$
of $X$\textit{\ can be extended to an }$\succcurlyeq $\textit{-increasing
and uniformly continuous map on }$X$\textit{.}

\bigskip

As every continuous map on a compact metric space is uniformly continuous,
an immediate consequence of Corollary 14 is the following observation which
provides a companion to Nachbin's extension theorem.

\bigskip

\noindent \textsf{\textbf{Corollary 15.}}\textbf{\ }\textit{Let }$%
(X,d,\succcurlyeq )$\textit{\ be a radial metric poset, and }$S$\textit{\ a
nonempty compact subset of }$X$\textit{. Then, for every }$\succcurlyeq $%
\textit{-increasing }$f\in C(S),$\textit{\ there is an }$\succcurlyeq $%
\textit{-increasing }$F\in UC(X)$ \textit{with }$F|_{S}=f$\textit{.}

\bigskip

At the cost of imposing the radiality property, this result drops the
topological requirement of being normally ordered in Nachbin's extension
theorem, and in addition, it guarantees the uniform continuity of the
extension as opposed to its mere continuity.

\subsection{The Monotone Uniform Extension Property}

It should be noted that the similarity of the statements of Theorem 4 and
Corollary 14 is misleading. To clarify this point, let us say that a
partially ordered metric space $(X,d,\succcurlyeq )$ has the \textit{%
monotone uniform extension property} if for every closed $S\subseteq X$ and $%
\succcurlyeq $-increasing and bounded $f\in UC(S),$ there is an $%
\succcurlyeq $-increasing $F\in UC(X)$ with $F|_{S}=f.$ The point we wish to
make is that this property is categorically different than the monotone
Lipschitz extension property. After all, the proof of Theorem 2 shows that a
finite metric poset has the monotone Lipschitz extension property iff that
metric poset is radial. In other words, finiteness of the carrier does not
allow us improve Theorem 2. By contrast, one can inductively prove that
every finite metric poset has the monotone uniform extension property. More
generally, an immediate application of Nachbin's extension theorem yields
the following fact: 

\medskip

\noindent \textit{Let }$(X,d,\succcurlyeq )$\textit{\ be a metric poset such
that }$(X,d)$\textit{\ is an UC-space.}\footnote{%
An \textit{UC-space} (also known as an \textit{Atsuji space}) is a metric
space such that every real-valued continuous function on it is uniformly
continuous. These spaces were originally considered by \cite{Atsuji}, \cite%
{M-P} and \cite{Nagata}, and were later studied extensively by \cite{Be2}, 
\cite{Be3} and \cite{J-K}, among others.) Various characterizations of
UC-spaces are known. For instance, a metric space $X=(X,d)$ is an UC-space
iff every open cover of it has a Lebesgue number, which holds iff $d(A,B)>0$
for every nonempty disjoint closed subsets $A$ and $B$ of $X.$}\textit{\
Then, }$(X,d,\succcurlyeq )$ \textit{has the monotone uniform extension
property, provided that it is normally ordered.}

\medskip

\noindent The family of all partially ordered metric spaces with the
monotone uniform extension property is thus much larger than that of radial
metric posets. Characterization of this family remains as an interesting
open problem.

\bigskip

\noindent \textsc{Acknowledgement }We thank professors Jerry Beer, Hiroki
Nishimura, Gil Riella, and Nik Weaver for their insightful comments at the
development stage of this work.

\end{document}